\documentclass{amsart}

\usepackage {amsfonts}
\usepackage[english]{babel}
\usepackage{color}
\def\G{\Gamma}
\def\d{\delta}
\def\a{\alpha}

\def\p{\varphi}
\def\e{\varepsilon}
\def\l{\lambda}
\def\L{\Lambda}

\def\t{\theta}

\def\s{\sigma}

\def\R{{\mathbb R}}

\def\N{{\mathbb N}}
\def\Z{{\mathbb Z}}

\DeclareMathOperator{\supp}{supp}
\DeclareMathOperator{\dist}{dist}
\DeclareMathOperator{\diam}{diam}

\newtheorem{Def}{Definition}
\newtheorem{Th}{Theorem}
\newtheorem{Pro}{Proposition}
\newtheorem*{Cor}{Corollary}

\begin{document}

\title{A new description of uniformly spread discrete sets}

\author{A.Dudko, S.Favorov}

\address{Sergii Favorov,
\newline\hphantom{iii}  V.N.Karazin Kharkiv National University
\newline\hphantom{iii} Svobody sq., 4, Kharkiv, Ukraine 61022}
\email{sfavorov@gmail.com}

\address{Artem Dudko,
\newline\hphantom{iii}  Institute of Mathematics of Polish  Academy of Sciences
\newline\hphantom{iii} Sniadeckich 8, 00-656, Warsaw, Poland}
\email{adudko@impan.pl}

\maketitle {\small
\begin{quote}
\centerline {\bf Abstract.}
We prove that each discrete set in the Euclidean space that has bounded changes  under every translation is a bounded perturbation of a square lattice, i.e., a uniformly spread set in the sense of Laszkovich.
In particular, the support of every Fourier quasicrystal with unit masses is uniformly spread.
\bigskip


\end{quote}
}

   \bigskip
\section{Introduction}\label{S1}
   \bigskip

One of the important parts of Laczkovich's proof of Tarski's famous problem on the equidecomposabilty of a square and a disk is the study of the so-called uniformly spread discrete sets in the plane.
Namely, M.Laszkovich in \cite[Theorem 3.1]{L1} proved that the following conditions are equivalent for any discrete set $A\subset\R^2$ and $\a>0$
\smallskip

\noindent
{\bf $a)$}  there exists a constant $C<\infty$ and a bijection $\s$ from $A$ onto $\a^{-1/2}Z^2$ such that $\sup_{x\in A} |\s(x)-x|<C$,
\vskip 0.1cm\noindent
{\bf $b)$} there exists a constant $C'<\infty$ such that $|\#(A\cap H)-\a m_2(H)|<C'm_1(\partial H)$ for every bounded Jordan domain $H$ with $\diam H\ge1$.
\vskip 0.1cm\noindent
Here, as usual, $\# E$ means the number of points of the finite set $E$, and $m_d(F)$ means the $d$-dimensional Lebesgue measure of the set $E$.

In  \cite{L2} M. Laszkovich showed that the direct analogue of this result is incorrect for $d>2$, and proved the equivalence  of the following conditions for a discrete set $A\subset\R^d,\ d>2$:
\smallskip

\noindent
$a'$)  there is a constant $C<\infty$ and a bijection $\s$ from $A$ onto $\a^{-1/d}Z^d$ such that $\sup_{x\in A} |\s(x)-x|<C$,
\smallskip

\noindent
$b'$) there is a constant $C'<\infty$ such that $|\#(A\cap H)-\a m_d(H)|<C'm_d\{x\in\R^d:\,\dist(x,\partial H)\le1\}$  for every bounded measurable set $H\subset\R^d$ with $\diam H\ge1$.
\smallskip

Here, a set is called discrete if its intersection with any ball is finite.

 M.Laszkovich called the sets  satisfying  conditions $a)$, $b)$ (or $a'$), $b'$) for $d>2$) {\it the uniformly spread sets}.

Notice that all conditions depend on the number $\a\in(0,\infty)$. It seems interesting to obtain a criterion for a set to be uniformly spread  that is independent of $\a$. This is what is done in this article.

The proof method implies a new statement about solutions of equations on infinite graphs, which may be useful in Graph Theory and its applications.

We also apply our result to Fourier quasicrystals with unit masses and show that their supports are uniformly spread.
\medskip

The paper is organized as follows. In Section \ref{S2} we give the necessary definitions and formulate the main results. Proofs of auxiliary Propositions \ref{P1}--\ref{P4} and Theorem \ref{T1} are given in the next Section \ref{S3}.
In Section \ref{S4} we present a proof of Theorem \ref{T2}. In Section \ref{S5} we generalize our results to multisets in $\R^d$ and apply  them to Fourier quasicrystals with unit or positive integer masses.
In the final section, we formulate questions related to roughly shift-invariant sets that seem interesting to us and for which we do not know the answers.

\bigskip
\section{Main definition and main results}\label{S2}
   \bigskip

The following  definition is the main one in our article:
\begin{Def}\label{D1}
A discrete set $A\subset \R^d$  is called a roughly shift-invariant set if there exist $L<\infty$ such that for any  $x\in\R^d$ there is a bijection $\s_x:\,A\to A$ such that
\begin{equation}\label{def}
  \sup_{a\in A}|a+x-\s_x(a)|<\infty.
\end{equation}
\end{Def}
Note that neither conditions $a)$ and $a'$) nor Definition \ref{D1}  depend on the choice of norm in $R^d$.
\medskip

Denote by $B(b,r),\ b=(b_1,\dots,b_d)\in\R^d$, the open ball of radius $r$ with center in $b$, and by $Q(b,r)$ the half-closed cube, i.e., the ball in $l^\infty$-norm of diameter $r$
with the vertex having minimal coordinates at $b$:
$$
Q(b,r)=\{x=(x_1,\dots,x_d)\in\R^d:\,b_j\le x_j<b_j+r,\ j=1,\dots,d\}.
$$
The set $E\subset\R^d$ is called {\sl relatively dense} if there exists $R<\infty$ such that $E\cap B(b,R)\neq\emptyset$ for every ball $B(b,R),\ b\in\R^d$.
Clearly, it suffices to verify \eqref{def} only for $x$ from any relatively dense set $E\subset\R^d$. The symbol $f(x)\sim g(x)$ as $x\to\infty$ means that  $\lim_{x\to\infty}f(x)/g(x)=1$.

Since $\a^{-1/d}\Z^d$ is roughly shift-invariant, we see that the conditions $a)$ or $a'$) for the set $A$ imply that $A$ is also roughly shift-invariant.

\begin{Th}\label{T1}
 For every roughly shift-invariant set $A$, there exists a density $D>0$ such that, uniformly with respect to $x\in\R^d$,
$$
  \lim_{T\to\infty}\frac{\#(A\cap B(x,T))}{ m_d(B(x,T))}=D.
$$
 \end{Th}

 \begin{Th}\label{T2}
 Every roughly shift-invariant set $A\subset\R^d$ is uniformly spread, and there exists a constant $C<\infty$ and a bijection $\Theta$ from $\N$ onto $D^{-1/d}Z^d$ such that
 $$
 \sup_{a\in A}|a-\Theta(a)|<C.
$$
\end{Th}
 Analogs of these results for measures were formulated and proved in \cite{DF}.
\medskip

The condition $b'$) implies  that for every uniformly spread set $A$ there exists a  constant $D>0$ (density of $A$) such that
\begin{equation}\label{D}
|\#(A\cap B(x,R))-Dm_d(B(x,R))|\le C'm_d(B(0,R+1\setminus B(0,R-1)),
\end{equation}
 hence, we obtain
\begin{Cor}
For every roughly shift-invariant set $A$ of density $D$, uniformly in $x\in\R^d$,
$$
\#(A\cap B(x,R))=Dm_d(B(x,R))+O(R^{d-1})\qquad (R\to\infty).
$$
\end{Cor}
Note that Theorems \ref{T1} and \ref{T2} were first formulated  in \cite{K}, but the proofs given there contain some gaps and inaccuracies, so we present new ones here.

\bigskip
\section{Auxiliary results and proof of Theorem \ref{T1}}\label{S3}
\bigskip

In what follows, we will assume that $A$ is a roughly shift-invariant set in $\R^d$ and there exists an integer $L$ such that
\begin{equation}\label{d}
\|a+x-\s_x(a)\|<L
\end{equation}
 for each $x\in\R^d$ and an appropriate bijection $\s_x:\,A\to A$.
Here $\|\cdot\|$ means the $l^\infty$ norm in $\R^d$.  Also, denote by $\textbf{1}$ the vector in $\Z^d$ with all coordinates equal to $1$.

\begin{Pro}\label{P1}
 There exists $K<\infty$ such that $\#(A\cap Q(x,1))<K$ for all $x\in\R^d$, and for $N\in\N$
$$
\#(A\cap Q(x,N))\le KN^d.
$$
If $N^{1/2}>d2^dLK$,  then
\begin{equation}\label{d1}
 \#(A\cap Q(x-L\textbf{1},N+2L)\setminus Q(x,N))<N^{d-1/2}.
\end{equation}
\end{Pro}
{\bf Proof}. For every $a\in A\cap Q(x,1)$ we have $a-x\in Q(0,1)$, hence, by \eqref{d}, we can associate a point  $a'\in A\cap Q(-L\textbf{1},2L+1)$ with $a$ such that $\|a'-(a-x)\|<L$, and distinct
points of $A\cap Q(-L\textbf{1},2L+1)$ correspond to different points of $A\cap Q(x,1)$. Hence, $\#(A\cap Q(x,1))\le \# A\cap Q(-L\textbf{1},2L+1):=K$. The second and third inequalities follow from the fact that
any cube $Q(x,N)$ is a union of $N^d$ pairwise disjoint cubes $Q(x^{(j)},1)$, and for $N$ large enough the set $Q(x-L\textbf{1},N+2L)\setminus Q(x,N)$ can be covered
 by $(N+2L)^d-N^d<2^ddL N^{d-1}<N^{d-1/2}/K$  cubes with unit edges.
\qquad$\Box$

\begin{Pro}\label{P2}
For $N$ sufficiently large, we have for each $x\in\R^d$
$$
  |\#(A\cap Q(x,N))-\#(A\cap Q(0,N))|<N^{d-1/2}.
$$
\end{Pro}

{\bf Proof}. Arguing as in the proof of Proposition \ref{P1}, we see that there is a one-to-one correspondence between $A\cap Q(x,N)$ and a subset of $A\cap Q(-L\textbf{1},N+2L)$.
Hence by \eqref{d1},
\begin{multline*}
  \#(A\cap Q(x,N))-\#(A\cap Q(0,N))\le\#(A\cap Q(-L\textbf{1},N+2L))-\#(A\cap Q(0,N))\\=\#(A\cap Q(-L\textbf{1},N+2L)\setminus Q(0,N))<N^{d-1/2}
\end{multline*}
for sufficiently large $N$.
In the same way, we estimate the difference $\#(A\cap Q(0,N))-\#(A\cap Q(x,N))$. \qquad $\Box$

\begin{Pro}\label{P3}
For  each $x\in\R^d$ the set $Q(x,2L)$ is nonempty, and for $N\in\N$ such that $N/(2L)$ is integer we get $\#(A\cap Q(x,N))\ge (N/(2L))^d$.
\end{Pro}
{\bf Proof}. Let  $a$ be an arbitrary point of $A$. Then $a+(x+L\textbf{1}-a)\in Q(x+L\textbf{1},L)$, therefore there exists $a'\in A\cap Q(x,2L)$.
It remains to note that every cube $Q(x,N)$ is a union of $(N/2L)^d$ cubes $Q(x^{(j)},2L)$.  \qquad $\Box$

\begin{Pro}\label{P4}
 There exists a number $D>0$ such that, uniformly with respect to $x\in\R^d$,
$$
  \lim_{T\to\infty} T^{-d}\#(A\cap Q(x,T))=D.
$$
 \end{Pro}

{\bf Proof}. Notice that any cube $Q(x,qN),\ x\in\R^d,\ N,\,q\in\N$, is a union of $q^d$ cubes $Q(x^{(j)},N)$ with some $x^{(j)}\in\R^d$, therefore by Proposition \ref{P2}
$$
     |q^{-d}\#(A\cap Q(x,qN))-\#(A\cap Q(x,N))|<N^{d-1/2}.
$$
Fix $\e>0$. Let $N_1,\ N_2\in\N$ be such that $N_1^{-1/2}<\e/4,\ N_2^{-1/2}<\e/4$.  The previous estimate implies:
$$
     \left|\frac{\#(A\cap Q(x,N_1 N_2))}{(N_1 N_2)^d}-\frac{\#(A\cap Q(x,N_1))}{N_1^d}\right|<N_1^{-1/2},
$$
$$
     \left|\frac{\#(A\cap Q(x,N_1 N_2))}{(N_1 N_2)^d}-\frac{\#(A\cap Q(x,N_2))}{N_2^d}\right|<N_2^{-1/2},
$$
therefore,
$$
     \left|\frac{\#(A\cap Q(x,N_1))}{N_1^d}-\frac{\#(A\cap Q(x,N_2))}{N_2^d}\right|<\e/2.
$$
It follows that $N^{-d}\#(A\cap Q(x,N))$, $N\in\mathbb N$, is a Cauchy sequence. Therefore, the limit $D(x)=\lim\limits_{N\to\infty}N^{-d}\#(A\cap Q(x,N))$ for $N\in\mathbb N$ exists.

 From Proposition \ref{P2} we have, for any $x,y\in\R^d$:
$$
    |N^{-d}\#(A\cap Q(x,N))-N^{-d}\#(A\cap Q(y,N))|\leq N^{-1/2}\to 0,
$$
when $N\to\infty$. The latter implies that $D(x)$ does not depend on $x$. Denote this common value of $D(x)$ by $D$. It follows from Proposition \ref{P3} that it is strictly positive.
Now, for arbitrary $T>1$, setting $N=[T]$ (the integer part of $T$) we have for $T\to\infty$
$$
\frac{\#(A\cap Q(x,N))}{N^{d}}\sim\frac{\#(A\cap Q(x,N))}{(N+1)^{d}}\leq\frac{\#(A\cap Q(x,T))}{T^{d}}\leq \frac{\#(A\cap Q(x,N+1))}{N^{d}}\sim\frac{\#(A\cap Q(x,N+1))}{(N+1)^{d}}.
$$
 Since both the leftmost and the rightmost terms of the above inequality converge to $D$, we obtain that $\lim\limits_{T\to\infty}T^{-d}\#(A\cap Q(x,T))=D$ (where now $T$ is arbitrary). \qquad $\Box$
\medskip

{\bf Proof of Theorem \ref{T1}}. Since $\diam Q(x,r)=r\sqrt{d}$, we can cover a ball $B(x,r^2)$ with the union of $M_r$ mutually disjoint cubes $Q(x_j,r)$  such that
$$
\bigcup_{j=1}^{M_r}Q(x_j,r)\subset B(x,r^2+\sqrt{d}r).
$$
 We have for $r\to\infty$
$$
m_d[B(x,r^2)]\sim{M_r}r^d\sim m_d[B(x,r^2+\sqrt{d}r)],
$$
and
$$
\#[A\cap B(x,r^2)]\le\sum_{j=1}^{M_r}\#[A\cap Q(x_j,r)]\le\#[A\cap B(x,r^2+\sqrt{d}r)].
$$
By proposition 4, the ratio $\#[A\cap Q(x,r)]/r^d$ tends to $D$ uniformly in $x\in\R^d$, hence
$$
\limsup_{r\to\infty}\frac{\#[A\cap B(x,r^2)]}{m_d(B(x,r^2))}\le\lim_{r\to\infty}\frac{1}{M_r}\sum_{j=1}^{M_r}\frac{\#[A\cap Q(x_j,r)]}{r^d}=D,
$$
and
$$
\liminf_{r\to\infty}\frac{\#[A\cap B(x,r^2+\sqrt{d}r)]}{m_d(B(x,r^2+\sqrt{d}r))}\ge\lim_{r\to\infty}\frac{1}{M_r}\sum_{j=1}^{M_r}\frac{\#[A\cap Q(x_j,r)]}{r^d}=D.
$$
The last two inequalities imply the assertion of the theorem. \qquad $\Box$

\bigskip
\section{Proof of Theorem \ref{T2}}\label{S4}
\bigskip

Let  $N\in\mathbb N$ be a number such that
\begin{equation}\label{L}
(N/(2L))^d>3^d(N^{d-1/2}+2),\qquad N^{1/2}>d2^dLK.
\end{equation}
For $i\in\Z^d$, set $P_i=\#(A\cap Q(iN,N))$. The idea of the proof of Theorem \ref{T2} is to construct (in several stages) a collection of integer numbers $t_{i,j},\  \|i-j\|=1$ such that
\begin{equation}\label{tij}
 N^d=P_i-\sum_{j:\|i-j\|=1}t_{i,j},\qquad t_{i,j}=-t_{j,i}\quad \forall\  i,j\in Z^d,
\end{equation}
 and
 \begin{equation}\label{ij}
  \sum\limits_{j:\|i-j\|=1}|t_{i,j}|\leq \min\{P_i, N^d\}\quad \forall\  i\in Z^d.
\end{equation}
Indeed, suppose that such a collection of integers exists. Take  a pair $(i',j')\in Z^d$ with $\|i'-j'\|=1$ and $t_{i',j'}\neq0$. Without loss of generality  suppose that $t_{i',j'}>0$.
Move $t_{i',j'}$ points of the set $A\cap Q(i'N,N)$ to the cube $Q(j'N,N)$. We may assume that each new point coincides with no point of $A$.
Denote the set $A\cap Q(j'N,N)$ with added $t_{i',j'}$ points as $A_{j'}$, and the set $A\cap Q(i'N,N)$ without these  points as $A_{i'}$. Set $P'_{j'}=\#(A_{j'}),\ P'_{i'}=\#(A_{i'})$. Since
$$
P'_{j'}=P_{j'}+t_{i',j'}=P_{j'}-t_{j',i'},\qquad P'_{i'}=P_{i'}-t_{i',j'},
$$
  we obtain
 $$
  N^d=P'_{i'}-\sum_{j\neq j':\|j-i'\|=1}t_{i',j},\qquad N^d=P'_{j'}-\sum_{i\neq i':\|i-j'\|=1}t_{j',i}.
 $$
 Replace the set $A\cap Q(j'N,N)$ with $A_{j'}$, and the set $A\cap Q(i'N,N)$ with $A_{i'}$. The set
 $$
  A'=A\setminus[A\cap(Q(j'N,N)\cup Q(i'N,N))]\cup A_{j'}\cup A_{i'}
 $$
 satisfies conditions \eqref{tij} with $P_{i'}$ replaced by $P'_{i'}$, $P_{j'}$ by $P'_{j'}$ and with the exclusion of $t_{i',j'}$ and $t_{j',i'}$.

 The inequalities \eqref{ij} show that we can repeat this procedure with all pairs $(i,j)$ such that $t_{i,j}\neq0$ and not move any point more than once. Therefore, every point of $A$ can be moved no more than $2N$
in the $l^\infty$ metric.  As a result, we obtain the set $\tilde A$ such that
$$
  N^d=\#[\tilde A\cap Q(iN,N)]\qquad \forall\ i\in\Z^d.
$$
Since each cube $Q(iN,N)$ contains exactly $N^d$ points of $\Z^d$, we see that $\tilde A$ satisfies $a')$ with $C=N$, and $A$ satisfies $a')$ with $C=3N$ in the metric $l^\infty$.
\medskip

{\bf First stage}.
 For $x\in\Z^d$, denote by $\s_x$ a bijection $A\to A$ such that $\|a+x-\s_x(a)\|<L$ for all $a\in A$.
Let $M^x_{j,i}$ be the set of points from  $Q(jN,N)$ that are preimages of points from $A\cap Q((i+x)N,N)$ under the bijection $\s_{Nx}$:
$$
M^x_{j,i}=Q(jN,N)\cap[\s^{-1}_{Nx}(A\cap Q((i+x)N,N))].
$$
 Since $\|a+Nx-\s_{Nx}(a)\|<L$ for $a\in Q(jN,N)$ and $N>2L$, we see that
\begin{equation}\label{d2}
\s_{Nx}(M^x_{j,i})\subset A\cap[Q((j+x)N-L\textbf{1},N+2L)\setminus Q((j+x)N,N)]\qquad\text{for}\quad\|j-i\|=1,
\end{equation}
and $\s_{Nx}(M^x_{j,i})=\emptyset$ for $\|i-j\|>1$.  Hence, $M^x_{j,i}=\emptyset$ for $\|i-j\|>1$.
Therefore $\s_{Nx}^{-1}$ gives a one-to-one correspondence between points of the set $A\cap Q((i+x)N,N)$ and  points of the set
\begin{equation}\label{M}
A\cap Q(iN,N)\cup\left[\bigcup_{j:\|i-j\|=1}M^x_{j,i}\right]\setminus\left[\bigcup_{j:\|i-j\|=1}M^x_{i,j}\right].
\end{equation}
Set $ p^x_{i,j}=\# M^x_{i,j}-\# M^x_{j,i}$. We get from \eqref{M} for all $i\in\Z^d$
\begin{equation}\label{p}
    P_{i+x}=P_i-\sum_{j:\|i-j\|=1}p^x_{i,j}.
\end{equation}
It follows from \eqref{d1} and \eqref{d2} that $\#M^x_{j,i}=\#\s_{Nx}(M^x_{j,i})<N^{d-1/2}$ for $\|i-j\|=1$, hence
\begin{equation}\label{p1}
    |p^x_{i,j}|\leq\max\{M^x_{i,j},\,M^x_{j,i}\}\le\begin{cases}N^{d-1/2}&\text{if }\|j-i\|=1,\\ 0 &\text{if }\|j-i\|\neq1.\end{cases}
\end{equation}
Further, it follows from \eqref{p} that for every $T\in\N$
\begin{equation}\label{t}
    \frac{\sum_{x\in \Z^d\cap Q(0,T)}P_{i+x}}{T^d}=P_i\frac{\sum_{x\in\Z^d,x\in Q(0,T)}}{T^d}-\frac{\sum_{x\in\Z^d,x\in Q(0,T)}\sum_{j:\|i-j\|=1}p^x_{i,j}}{T^d},
\end{equation}
The ratio
$$
    T^{-d} \sum_{x\in\Z^d, x\in Q(0,T)}p^{x}_{i,j}
$$
is uniformly bounded for $T\in\N$.  Using the diagonal process, we can find a sequence $T_n\to\infty$ such that for all $i,\,j$
the ratios
$$
    T_n^{-d} \sum_{x\in\Z^d,x\in Q(0,T_n)}p^{x}_{i,j}
$$
have limits $p_{i,j}$ as $n\to\infty$. Note that $p_{i,j}=-p_{j,i}$ for all $i,\,j\in\Z^d$.

Next, we have
$$
 \#(A\cap Q(iN,TN))=\sum_{x\in \Z^d\cap Q(0,T)}\#(A\cap Q((i+x)N,N)).
$$
Therefore by Proposition \ref{P4} and our assumption
 $$
   \lim_{T\to\infty}\frac{\sum_{x\in \Z^d\cap Q(0,T)}P_{i+x}}{T^d}=\lim_{T\to\infty}\frac{N^d\#(A\cap Q(iN,TN))}{(TN)^d}=N^d.
 $$
 Since $\#[Q(0,T)\cap \Z^d]=T^d$, we obtain from \eqref{t}
 \begin{equation}\label{m}
   N^d=P_i-\sum_{i:\|i-j\|=1}p_{i,j}.
\end{equation}
Taking into account that $\#\{j:\|j-i\|=1\}<3^d$ for all $i\in\Z^d$, we get from \eqref{p1}
$$
\sum_{j:\|i-j\|=1}|p^x_{i,j}|<3^d N^{d-1/2}\qquad \forall\ x,\,i\in\Z^d,
$$
and
\begin{equation}\label{p2}
\sum_{j:\|i-j\|=1}|p_{i,j}|\le 3^d N^{d-1/2}\qquad \forall\ i\in\Z^d.
\end{equation}
We see that the numbers $p_{i,j}$ would be natural candidates for the numbers $t_{i,j}$ (see \eqref{tij}),  except that $p_{i,j}$ are not integer. In subsequent stages we will modify $p_{i,j}$ accordingly.
\medskip

{\bf Second stage}. Denote by $[y]$ the integer part and by $\{y\}$ the fractional  part of $y\in\R$. Initially, set $q_{i,j}=p_{i,j}$ for each $i,j\in Z^d$.
Let us say that a finite sequence of pairwise distinct points $i_0,i_1,\dots,i_{n-1}\in\Z^d$ forms a bad cycle if  $\|i_k-i_{k+1}\|=1$ and $q_{i_k,i_{k+1}}\not\in\Z$ for all $k$, where we set $i_n=i_0$.
Observe that there are at most a countable number of bad cycles. We will go through each bad cycle one by one and perform the following operation.

Let $\t=\min\{\{p_{k,k+1}\}:0\leq k<n\}$ for some bad cycle $i_0,i_1,\dots,i_{n-1}$. Clearly, $0<\t<1$. Replace all numbers $q_{i_k,i_{k+1}}$ with
$q_{i_k,i_{k+1}}=p_{i_k,i_{k+1}}-\t$, and the numbers $q_{i_{k+1},i_k}$ with $q_{i_{k+1},i_k}=p_{i_{k+1},i_k}+\t$. We get, $q_{i_{k+1},i_k}=-q_{i_k,i_{k+1}}$ for all $k$.
Therefore, if for some $i\in\Z^d$ we have $i=i_k$, then among $q_{i,j}$ with $\|i-j\|=1$ only  $q_{i_k,i_{k+1}}$ and $q_{i_k,i_{k-1}}$  are changed and
$$
   q_{i_k,i_{k+1}}+q_{i_k,i_{k-1}}= (p_{i_k,i_{k+1}}-\t)+(p_{i_k,i_{k-1}}+\t)= p_{i_k,i_{k+1}}+p_{i_k,i_{k-1}}.
$$
 Thus, after such an operation, we still have
\begin{equation}\label{q1}
 N^d=P_i-\sum_{j:\|i-j\|=1}q_{i,j},\qquad q_{j,i}=-q_{i,j},\qquad \forall i,\,j\in\Z^d.
\end{equation}

Observe that the above operation eliminates at least one bad cycle and does not create any new ones. Repeat this operation with each remaining bad cycle.
After a finite or infinite number of such operations, we find that there are no more bad cycles, and the equality \eqref{q1} remains valid. Moreover, each number $q_{i,j}$ is changed  by no more than $1$, hence
\begin{equation}\label{q2}
\sum_{j:\|i-j\|=1}|q_{i,j}|\le \sum_{j:\|i-j\|=1}|p_{i,j}|+3^d\qquad \forall i\in\Z^d.
\end{equation}

{\bf Third stage}.  Consider the graph $\G$ with vertices $i\in\Z^d$ and edges $(i,j)$ such that $q_{i,j}\not\in\Z$.  Let $\G'$ be its connected component that contains more than one vertex.
Enumerate vertices of $\G'$ by numbers $0,1,2,\dots$, such that $i_0$ is arbitrary and for every $k>0$ there is a vertex $i_l$ with $l<k$ such that the edge $(i_l,i_k)\in\G'$.
According to the second stage, the graph $\G$ has no cycles; hence, for every $k$ there is only one index $l$ with this property.

Replace $q_{i,j}$ with $t_{i,j}$ for each pair $(i,j)\in Z^d\setminus\G$. We will change the remaining $q_{i,j}$ to integers inductively. Suppose that for some $k$ we have replaced all numbers $q_{i_l,i_m},\ l\le k,\,m\le k$, with integers  $t_{i_l,i_m}$ such that for every $l\le k$
\begin{equation}\label{q3}
                                -1+\left(P_{i_l}-N^d-\sum_{j\in\Z^d\setminus\G'}t_{i_l,j}\right)<\sum_{m\le k}t_{i_l,i_m}+\sum_{m>k}q_{i_l,i_m}<\left(P_{i_l}-N^d-\sum_{j\in\Z^d\setminus\G'}t_{i_l,j}\right)+1.
\end{equation}
On the other hand, we have not changed any $q_{i_l,i_m}$ for $l>k$, so by \eqref{q1} for such $l$
\begin{equation}\label{q4}
\sum_{m>k} q_{i_l,i_m}=P_{i_l}-N^d-\sum_{j\in\Z^d\setminus\G'}t_{i_l,j}.
\end{equation}
As we noted above, for the vertex $i_{k+1}$ there exists a unique index $l<k$ such that $i_l\in\G'$ and $\|i_l-i_{k+1}\|=1$. Clearly, we can replace $q_{i_l,i_{k+1}}$ with
$t_{i_l,i_{k+1}}=[q_{i_l,i_{k+1}}]$ or $[q_{i_l,i_{k+1}}]+1$ in \eqref{q3} so that the inequality \eqref{q3} would turn into
$$
                              -1+\left(P_{i_l}-N^d-\sum_{j\in\Z^d\setminus\G'}t_{i_l,j}\right)<\sum_{m\le k}t_{i_l,i_m}+t_{i_l,i_{k+1}}+\sum_{m>k+1}q_{i_l,i_m}<\left(P_{i_l}-N^d-\sum_{j\in\Z^d\setminus\G'}t_{i_l,j}\right)+1,
$$
which is equivalent to \eqref{q3} with $k+1$ instead of $k$. Next, replace $q_{i_{k+1},i_l}$ with $t_{i_{k+1},i_l}=-t_{i_l,i_{k+1}}$. Since
$$
|t_{i_{k+1},i_l}-q_{i_{k+1},i_l}|=|t_{i_l,i_{k+1}}-q_{i_l,i_{k+1}}|<1,
$$
we see that the equality \eqref{q4} implies inequality \eqref{q3} for  $k+1$ instead of $k$ and with $l=k$.
These reasoning is valid also for $k=0$, since in this case \eqref{q3}  has the form
$$
 P_{i_0}-N^d-\sum_{j\in\Z^d\setminus\G'}t_{i_0,j}-1<\sum_{m>0}q_{i_0,i_m}<P_{i_0}-N^d-\sum_{j\in\Z^d\setminus\G'}t_{i_0,j}+1,
$$
which is trivially true in view of \eqref{q1} with $i=i_0$. By an induction argument, we can replace all non-integer numbers $q_{i,j},\ (i,j)\in\G'$ with integers $t_{i,j}$ such that \eqref{q3} satisfies for all $k$.

Further, for every fixed $l$, there is $k$ such that $\|i_l-i_m\|>1$ for all $m>k$. In this case the inequality \eqref{q3} turns into
$$
 P_{i_l}-N^d-\sum_{j\in\Z^d\setminus\G'}t_{i_l,j}-1<\sum_{m\le k}t_{i_l,i_m}<P_{i_l}-N^d-\sum_{j\in\Z^d\setminus\G'}t_{i_l,j}+1.
$$
All numbers in this inequality are integers, hence, we obtain \eqref{tij} for $i=i_l$ and consequently for all $i\in\G'$. Repeating this procedure for other connected components of $\G$,
we obtain  \eqref{tij} for all $i,j$.

 Furthermore, we change only non-integer $q_{i,j}$, and at most by 1, hence, it follows from Proposition \ref{P3}, \eqref{L}, \eqref{p2}, and \eqref{q2}  that for all $i\in\Z^d$
\begin{equation}\label{t2}
\sum_{j:\|i-j\|=1}|t_{i,j}|\le \sum_{j:\|i-j\|=1}|q_{i,j}|+3^d\leq 3^d(N^{d-1/2}+2)<\min\{P_i,N^d\}.
\end{equation}
The argument at the beginning of the proof of Theorem \ref{T2} completes the proof. \qquad$\Box$
\medskip

Notice that, following the proof of Theorem \ref{T2}, we can obtain a new statement about solutions of equations on general graphs.
For a directed graph $\Gamma=(V,E)$, denote by $E_v$ the set of edges starting at $v$. For a directed edge $e$, denote by $-e$ the same edge taken in the opposite direction.
Similarly to stages 2 and 3 of the proof we can show the following:

\begin{Th}\label{T3}
 Let $\Gamma=(V,E)$ be any directed graph with a countable number of vertices and finite degree at each vertex. Let $R_v$, $v\in V$, be a collection of integer numbers.
  Assume that there exist real numbers $p_e,\,e\in E$, such that  for every $e\in E$ one has $p_e=-p_{-e}$, and for every $v$ one has $\sum\limits_{e\in E_v}p_e=R_v$.
  Then there exit integer numbers $t_e\in \{[p_e],[p_e]+1\}$ such that
  the same identities hold:
 $$
 \sum\limits_{e\in E_v}t_e=R_v\;\;\text{for every}\;\;v\in V,\qquad t_e=-t_{-e}\;\;\text{for every}\;\;e\in E.
 $$
 \end{Th}

\bigskip
\section{Roughly shift-invariant multisets and applications to Fourier quasicrystals}\label{S5}
   \bigskip

We may extend the concept of discrete sets to discrete multisets, which are discrete sets such that a natural number (multiplicity) is assigned to each point of the set. It is more convenient to describe such sets
as a sequence   $A=\{a_n\}_{n\in\N}\subset\R^d$ that has no finite limit points. To extend the conditions $a)$ and $a')$ to multisets, we should replace the condition
$$
\sup_{x\in A} |\s(x)-x|<C\quad \text{for some bijection}\ \s:\,A\to\a^{-1/d}Z^d
$$
with
$$
 \sup_{n\in\N}|\s(n)-a_n|<C\quad \text{for some bijection}\ \s:\,\N\to\a^{-1/d}\Z^d.
$$
The definition of a roughly shift-invariant multiset has the following form:
\begin{Def}\label{D2}
A discrete multiset $A=\{a_n\}_{n\in\N}\subset \R^d$  is a roughly shift-invariant set if  for any  $x\in\R^d$ there is a bijection $\s_x:\,\N\to\N$ such that
$$
  \sup_{n\in\N}|a_n+x-a_{\s_x(n)}|<\infty.
$$
\end{Def}
Now, let the symbol $\# E$ denote the number of points of a finite multiset $E$ where each point is counted according to its multiplicity. Our theorems take the following form for multisets:

\begin{Th}\label{T4}
 For every roughly shift-invariant multiset $A$, there exists a density $D>0$ such that uniformly with respect to $x\in\R^d$
$$
  \lim_{T\to\infty}\frac{\#(A\cap B(x,T))}{ m_d(B(x,T))}=D.
$$
 \end{Th}

 \begin{Th}\label{T5}
 Every roughly shift-invariant multiset $A\subset\R^d$ is  uniformly spread, and there is a constant $C<\infty$ and a bijection $\Theta$ from $\N$ onto $D^{-1/d}Z^d$ such that
 $$
 \sup_{n\in\N}|a_n-\Theta(n)|<C.
$$
\end{Th}
Indeed, take a sequence $r_n\in\R^d$ such that $|r_n|<1$ and $a_n+r_n\neq a_k+r_k$ for each $n\neq k$.
 Then, the map  $\t:\, n\to a_n+r_n$ is a bijection from $\N$ to the set $A'=\{a_n+r_n\}_{n\in\N}$ such that $|a_n-\t(n)|<1$,
and Theorems \ref{T4} and \ref{T5} follow immediately from Theorems \ref{T1} and \ref{T2}. To verify the Corollary for $n>1$, one can choose points $r_n$ such that $|a_n+r_n|=|a_n|$.
For $n=1$, one can use Proposition \ref{P1}, which is also valid for multisets.
\bigskip

The notion of a roughly shift-invariant (multi)set is an extension of the notion of an almost periodic set in $\R^d$:
\begin{Def} [\cite{L} Appendix VI for $d=1$, \cite{FK} for $d>1$]\label{D3}
A discrete multiset $A=\{a_n\}_{n\in\N}\subset\R^d$ is  almost periodic if for any $\e>0$ there is
a relatively dense set $E_\e$ of {\sl $\e$-almost periods} of $A$
$$
  E_\e=\{x\in\R^d:\,\sup_n|a_n+x-a_{\s_x(n)}|<\e\quad \text{for some bijection}\ \s_x:\,\N\to\N\}.
$$
\end{Def}
It is easy to see that every almost periodic multiset in $\R^d$ is a roughly shift-invariant multiset.
\medskip

  Almost periodic sets are naturally related to almost periodic measures and Fourier quasicrystals.
Recall that a complex measure $\mu$ on the Euclidean space $\R^d$ with discrete support (i.e., its intersection with any compact set is finite)
 is called a {\it Fourier quasicrystal} if $\mu$ is a tempered distribution, its Fourier transform in the sense of distributions $\hat\mu$
 is also a measure with discrete support, and both measures $|\mu|$ and $|\hat\mu|$ are tempered distributions as well.
 Here, we denote by $|\nu|(E)$ the variation of the complex measure $\nu$ on the set $E$.
 In fact, each Fourier quasicrystal is a form of some Poisson formula. Both of these objects are used and studied very actively, see, for example \cite{AKV}, \cite{LT}, \cite{OU}.
 In particular, Poisson formulas  were used by D. Radchenko and M. Viazovska in \cite{RV}.

Of greatest interest are Fourier quasicrystals with unit masses
\begin{equation}\label{a}
\mu=\sum_{\l\in\L}\d_\l,\quad \l\in\R^d,
\end{equation}
 ($\d_\l$, as usual, means the unit mass at the point $\l$). First nontrivial example of such Fourier quasicrystal was given by P.Kurasov and P.Sarnak \cite{KS}.
 Then A. Olevsky and A. Ulanovsky in \cite{OU} proved that for $d=1$ the support of any Fourier quasicrystal  \eqref{a}
is the zeros of real-root exponential polynomial, and vice versa, the zeros of any real-root exponential polynomial is the support of such Fourier quasicrystal.
W. Lawton and A. Tsikh \cite{LT} partially extended this result to Euclidean spaces  $\R^d,\ d>1$.

It was proved in \cite {M}, \cite{F1} that the convolution $\mu\star\p$ of every non-negative Fourier quasicrystal $\mu$  with any compactly supported continuous function $\p$
is an almost periodic function.  By definition, this means that the measure $\mu$ is almost periodic. If $\mu$ has only unit masses, then $\supp\mu$ is an almost periodic set
(for $d=1$ see \cite{FRR}, for $d>1$ see \cite{FK}). On the other hand, Theorem \ref{T2}  yields that each  almost periodic set in $\R^d$ is  uniformly spread.
Using  the Corollary, we obtain the following result:

\begin{Th}\label{T6}
For the support $\L$ of any Fourier quasicrystal \eqref{a} in $\R^d$   we have, uniformly in $x\in\R^d$,
$$
\#(\L\cap B(x,R))=Dm_d(B(x,R))+O(R^{d-1})\qquad(R\to\infty).
$$
Moreover, there is  a bijection  $\s:\,D^{-1/d}Z^d\to \L$ such that for every $\l=\s(h)\in\L$, uniformly in $h\in D^{-1/d}Z^d$, one has:
$$
 \l=h+O(1)
$$
\end{Th}
The first equality  was obtained  earlier by other methods in \cite{AKV}.
\medskip

If a Fourier quasicrystal $\mu$ has positive integer masses, we assume that the multiplicity of each point $\l\in\supp\mu$ is equal to the mass $\mu\{\l\}$,
then $\supp\mu$ becomes an almost periodic multiset. Using  Theorem \ref{T5} instead of Theorem \ref{T2}, we get that each  almost periodic multiset in $\R^d$ is uniformly spread.
Therefore in this case the both statements of Theorem \ref{T6} are also valid.

Finally, let a Fourier quasicrystal $\mu=\sum_{\l\in\L}a_\l\d_\l$ have uniformly discrete support $\L$ and complex masses $a_\l$ such that $\inf_\l|a_\l|>0$. By \cite[Lemma 4 and Theorem 11]{F0},
the measure $\mu$ is almost periodic. By \cite[Lemma 3]{F0}, the set $\L$ is almost periodic too,  therefore it satisfies both statements of Theorem \ref{T6} as well.

\bigskip
\section{Some questions}\label{S6}
   \bigskip

{\bf Question 1}. It was proved in \cite{F2} that every almost periodic multiset $A=\{a_n\}_{n\in\Z}\subset\R$ under appropriate numbering has the type
$$
a_n=n/D+\phi(n)\quad\text{with an almost periodic mapping}\quad \phi:\,\Z\to\R.
$$
Is there an analog of this result for $\R^d,\ d>1$?
\medskip

{\bf Question 2}.  What is the optimal value of the perturbation constant $C$?

\end{document}